\documentclass{article}

\usepackage[final]{neurips_2020_tda_style}
\def\sfb{\sffamily \bfseries}
\usepackage{booktabs}

\usepackage[utf8]{inputenc} 
\usepackage[T1]{fontenc}    
\usepackage{amsfonts}       
\usepackage{booktabs}       
\usepackage{hyperref}       
\usepackage{url}            
\usepackage{microtype}      
\usepackage{nicefrac}       
\usepackage{paralist}       
\usepackage{siunitx}        
\usepackage{graphicx}
\title{Hotspot identification for Mapper graphs}

\author{
  Ciara F. Loughrey$^{\text{\sfb 1}}$, Nick Orr$^{\text{\sfb 2}}$, Anna Jurek-Loughrey$^{\text{\sfb 1}}$ \\
  $^{\text{\sf 1}}$School of Electronics, Electrical Engineering and Computer Science\\
  $^{\text{\sf 2}}$Patrick G Johnston Centre for Cancer Research\\
  Queen's University Belfast, Northern Ireland\\
  \texttt{ \{cloughrey11, nick.orr, a.jurek\} @qub.ac.uk} \\
  \And
  Pawe\l{} D\l{}otko \\
  Dioscuri Centre in Topological Data Analysis \\
  Mathematical Institute, Polish Academy of Sciences
  Warsaw, Poland \\
  \texttt{pdlotko@impan.pl} \\
}

\begin{document}

\maketitle

\begin{abstract}
 Mapper algorithm can be used to build graph-based representations of high-dimensional data capturing structurally interesting features such as loops, flares or clusters. The graph can be further annotated with additional colouring of vertices allowing location of regions of special interest. For instance, in many applications, such as precision medicine, Mapper graph has been used to identify unknown compactly localized subareas within the dataset demonstrating unique or unusual behaviours. This task, performed so far by a researcher, can be automatized using hotspot analysis. In this work we propose a new algorithm for detecting hotspots in Mapper graphs. It allows automatizing of the hotspot detection process. We demonstrate the performance of the algorithm on a number of artificial and real world datasets. We further demonstrate how our algorithm can be used for the automatic selection of the Mapper lens functions. 
\end{abstract}

\section{Introduction}
The aim of hotspot analysis is to detect anomalous regions, such as subgroups of patients with distinctive survival, within a dataset. Those regions could indicate a novel phenomena (e.g. the existence of an unknown disease type (\citet{Nicolau2011, Lum2013, Cho2019}). The Mapper algorithm (summarized in~\ref{sec:mapper}) generates a graph representation of the layout of multidimensional data samples. The original data points are mapped to the low-dimensional space according to the \emph{lens} function of the algorithm, also known as the filter function. Once an appropriate lens function is used, the hotspots in the data and hotspots in the Mapper graph will be in correspondence. Consequently Mapper graphs provide a handy way to detect hotspots. There are however two issues; Firstly, manual analysis of large graphs, that are often obtained for large datasets, is prohibitive. Secondly, selection of appropriate parameters for Mapper construction that reveal the hotspots in data is nontrivial and typically requires the construction of multiple Mapper graphs that need to be analyzed. In order to address this challenge we propose a new technique for automatic detection of hotspots in graphs. Using a real-world breast cancer dataset, we demonstrate how the proposed algorithm can be used to automatically select a lens function based on its ability to discriminate subgroups of patients that present increased survival outcomes.

\section{Problem Statement} 
\label{PS}
In this work we address the problem of detecting regions within a Mapper graph that are structurally coherent and homogeneous on a value attribute of interest (e.g. survival), while also differing sufficiently from its neighbourhood within graph. A \emph{pullback} of such a hotspot region of a Mapper graph will then indicate a hotspot in the initial dataset. Let $X=\{x_1, \dots, x_k\}, x_i \in \mathbb{R}^n$ be the considered point cloud and $A : X \rightarrow \mathbb{R}$ indicate an attribute of data (e.g. survival of patients). Let $G={<}V,E{>}$ be a Mapper graph obtained from $X$ using chosen parameters of the constructions and $\hat{A} : V \rightarrow \mathbb{R}$ be the induced attribute function (as described in supplementary material). $C \subset V$ will be called a \emph{hotspot} within $G$ with respect to $\hat{A}$ if the following conditions holds:
\begin{enumerate}
	\item \textit{Connectedness:} any two vertices $\{v_i,v_j\} \in C$ are connected by a sequence of edges (path) from $E$ supported in $C$.
	\item \textit{Internal Homogeneity:} the dispersion of $\hat{A}$ values for data points across all vertices from $C$ is not more that a predefined threshold $\tau$. Formally, $\forall_{v_1,v_2 \in C} |\hat{A}(v_1) - \hat{A}(v_2)| \leq \tau$. 
	\item \textit{Neighbourhood Heterogeneity:} the values of $\hat{A}$ on $C$ are sufficiently different from the values of $\hat{A}$ on vertices within the neighbourhood of $C$ (denoted as $N_C$). We define $N_C$ as:
	\begin{equation} \label{eqn}
	N_C=\{C_i \text{ connected } |C_i\subset V: C_i \cap C = \emptyset, \exists_{v\in C}, \exists_{v'\in C_i}, {<}v,v'{>}\in E \}
	\end{equation}
	Then we require that $|\hat{A}(C) - \hat{A}(N_{C})| > \epsilon$. This definition is explained further in Figure~\ref{fig:neighbourhood_sup}.
	\item \emph{Size:} the size of $C$ (denoted as $S(C)$) is the size of a vertex, i.e. the number of points that are covered by the corresponding cluster. $S(C)$ should be large enough so that $C$ is not an outlier, but it is proportionally smaller than the size of its neighbourhood. Formally, $S(C) > \sigma_1$ and $S(N_{C}) -S(C)> \sigma_2$. Depending on the dataset, we may prefer to consider $S(C)$ as the total number of data points across all vertices in $C$ or the total number of nodes within $C_i$, or we can use both as the criteria.
\end{enumerate}




\section{Hotspot Detection Algorithm}
\label{sec:hotspot_detection}
The proposed algorithm proceeds in two steps; Step 1, non-intersecting connected components of $G$ that are homogeneous with respect to $\hat{A}$ are chosen. Note that there are cases when those regions are not uniquely defined, as described in ~\ref{sec:non_uniqueness_interlal_homogenity}. Therefore we chose any subset that satisfies the internal homogeneity criteria. Step 2, each component is classified as either hotspot or non-hotspot. 

Cluster detection can be explained as follows; Given $G$, the vertices are assigned an average value of an attribute function $\hat{A}$ (e.g. survival) – the average is computed over the points that are covered by a given vertex. To separate the regions with vastly different values we construct a function on edges, $F'$, that captures a large gradient of $\hat{A}$ over the edges of the graph. We then build a dendrogram based on the $F'$ of edges. Put simply, two vertices of $G$ connected by an edge with similar $\hat{A}$ values will be merged quickly, as the filtration of the edge will be small. When we are given two vertices joined by an edge with very different values of $\hat{A}$, the edge joining them will appear late in the filtration. The process of building dendrograms requires the filtration on vertices to be 0 for the construction to make sense.

We assume that $\tau$, $\epsilon$ and $\sigma_1$ should be set by the user as they strongly depend on the domain and $F'$. In our experiments we manually set values of $\epsilon$ and $\sigma_1$ while the optimal value of $\tau$ is determined from the dendrogram.

\textbf{Step 1: Cluster detection} (Assuring the connectedness of $G$ and the internal homogeneity conditions)
\begin{enumerate}
  \item Define a new function $F'$ on $V$ and $E$: 
	\begin{equation}
  	\forall_{v\in V}\, F'(v)=0, \quad \forall_{{<}v_{i},v_{j}{>}\in E}\, F'(<v_{i},v_{j}>)=|\hat{A}(v_{i})-\hat{A}(v_{j})|
  	\label{eq:three}
	\end{equation}
  \item Perform single linkage on $V$ using $F'$ and obtain the corresponding dendrogram. 
  \item Identify all connected components $C_1, \dots , C_n$ that are connected in the dendrogram below the level $\tau$. Set $\tilde{C}=\{C_1, \dots , C_n\}$.
\end{enumerate}

\textbf{Step 2: Cluster classification} (Assuring the size and the neighbourhood heterogeneity conditions)
\begin{enumerate}
	\item For each $C_i \in \tilde{C}$, calculate size of $N_{C_i}$ as:
		\begin{equation}
	 	S(N_{C_i})=\frac{1}{|N_{C_i}|}\sum_{C\in N_{C_i}}S(C)
	 	\label{neighbourhood_eq}
		\end{equation}	
	\item If $S(C_i) < \sigma_1$ or $|S(N_{C_i}) -S(C_i)|< \sigma_2$ then classify $C_i$ as a non-hotspot and remove it from $\tilde{C}$. We assume that if $C_i$ is very small then it should be considered as an outlier rather than a hotspot. At the same time for $C_i$ to be a hotspot it should be proportionally smaller than  its neighbourhoods. We propose for $\sigma_1$ to be a parameter set by the user, and for $\sigma_2$ to be calculated as one median absolute deviation of $\{S(C_1)\dots S(C_n)\}$.   
	\item For each $C_i$ calculate $\hat{A}(N_{C_i})$ as the mean value of $\hat{A}$ across all vertices within $N_{C_i}$
	\item If $|\hat{A}(C_i) - \hat{A}(N_{C_i})| > \epsilon$ then $C_i$ is considered as a hotspot.
\end{enumerate}

\section{Experimental Evaluation}
The hotspot detection algorithm was evaluated with Mapper graphs constructed with one toy dataset (Figure \ref{fig:artificialdatasets}a) and one real-word breast cancer dataset obtained from The Cancer Genome Atlas \citep{TCGA} as well as two large and complex artificial graphs (Figure~\ref{fig:artificialdatasets}b and~\ref{fig:artificialdatasets}c).


\subsection{Artificial Datasets} 
\label{Artifical section}
\begin{figure}[h]
\centering
\includegraphics[width=1\textwidth]{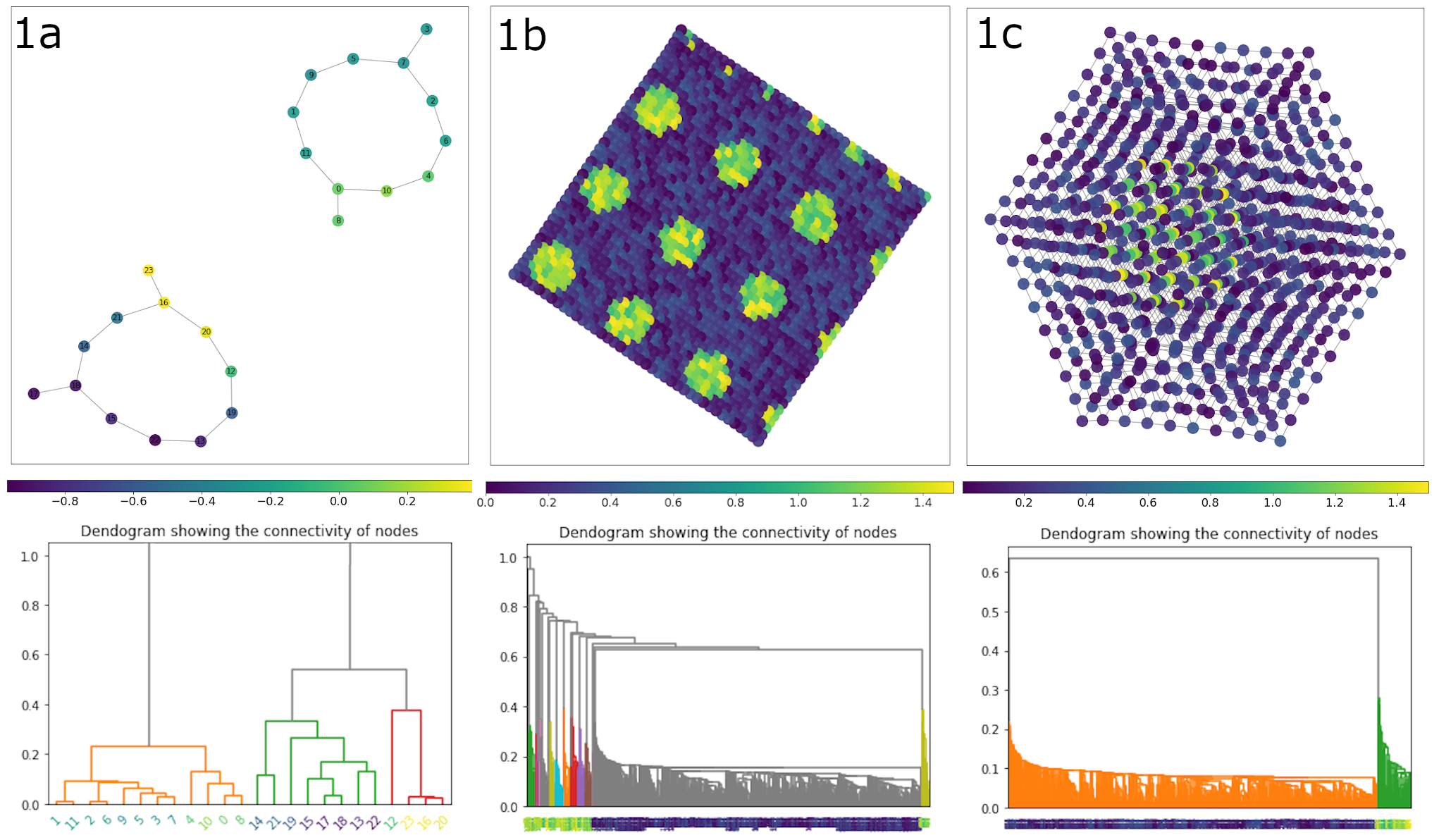} 
\caption{Three Mapper graphs used in the experiments and the corresponding dendrograms constructed with the function $F'$ from Equation \ref{eq:three}. In 1a, the node colour is based on the minimum value for vector contained in the node. In 1b and 1c, the method for determining the node colour is described in Section~\ref{Artifical section}}
\label{fig:artificialdatasets}
\end{figure}
The first dataset (Figure \ref{fig:artificialdatasets}a) consists of two 2-D circles constructed using the python “scikit-learn” package \citep{scikit-learn}. We artificially assigned values to the data points to obtain a hotspot on one of the circles, by defining $\hat{A}$ as the minimum value of each vector. A Mapper graph was built using L2-norm as the lens function with the number of bins and the overlap parameters set to $7$ and $20$, respectively. Clustering was performed using wards linkage with $6$ clusters per interval. Following phase one of the algorithm, three distinct candidates were identified. We set $\sigma_1$ = $2$ for nodes and the $\epsilon$ threshold at $0.01$. A single connected component corresponding to one circle was split into two regions, and the smaller yellow region of high values was identified as a hotspot (Figure \ref{fig:artificialdatasets}a). The two artificial graphs were obtained in the following way. Firstly, a sufficiently dense 2-D (Figure \ref{fig:artificialdatasets}b) or 3-D (Figure \ref{fig:artificialdatasets}c) grid of points is selected. Then all neighbouring grid elements are connected. In addition, a smaller number of connections between random vertices is added. Subsequently the function $\hat{A}$ is defined on a graph being a small random variable plus a correction. In the 2-D instance a correction of $1$ is added to all the grid points $(x,y)$ for which $\sin(x) + \sin(y) \geq 1$. In the 3-D case the correction is added to all grid points $(x,y,z)$ for which $x^2+y^2+z^z \leq 1$. No other corrections are added.

In the graph with multiple hotspots visualized in Figure \ref{fig:artificialdatasets}b, $9$ hotspots were identified according to a threshold of $\sigma_1$ = $15$ for nodes and $\epsilon$ = $0.01$. We observed that the algorithm was sensitive with respect to the specified minimum size of the hotspot. Setting the parameter too low resulted in larger numbers of very small hotspots that may be considered as outliers. The single hotspot in the 3-D graph in Figure \ref{fig:artificialdatasets}c was detected based on a minimum $\sigma_1$ value of $10$ for nodes and $\epsilon$ value of $0.01$. The results for each artificial dataset are summarised in Table~\ref{artificial_table}. Within the artificial graphs, no false positives or false negatives were detected by the hotspot analysis. The minimum lens function difference required some exploration to ensure a stringent threshold boundary at which to consider a candidate a hotspot.

\begin{table}
  \caption{%
   Results from hotspot detection on the artificial datasets. We describe the number of hotspots found for each dataset and the parameter settings for the hotspot detection algorithm. 
  }
  \label{artificial_table}
  \centering
  \begin{tabular}{lrrr}
    \toprule
    \textbf{Dataset} & \textbf{$\sigma_1$ (nodes)} & \textbf{$\epsilon$} & \textbf{Hotspot count} \\
    \midrule
    \textbf{Two circles}      & 2                   & 0.1                  & 1                                    \\
    \textbf{2-D graph}        & 15                  & 0.01                 & 9                                    \\
    \textbf{3-D graph}        & 10                  & 0.01                 & 1                                  \\
    \bottomrule
  \end{tabular}
\end{table}

\subsection{Real World Datasets}

\begin{figure}[h!]
\centering
\includegraphics[width=0.7\textwidth]{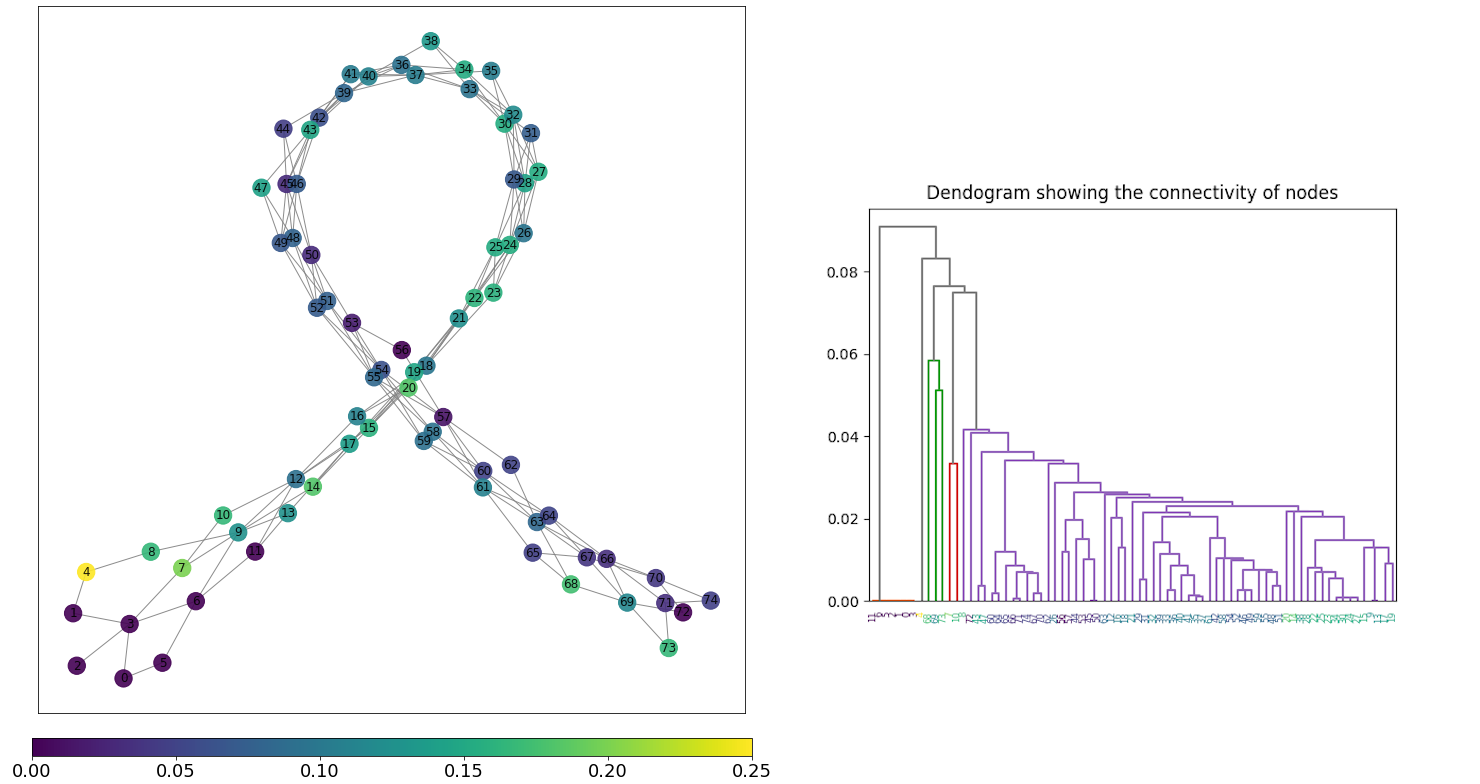} 
\caption{Mapper graph constructed from the TCGA dataset using the proposed algorithm (constructed with linear combination of a subset of features as the lens function and coloured by survival outcome), with the corresponding dendrogram of node connectivity based on edge weights}
\label{fig:linearsub}
\linespread{0.5}
\end{figure}
The real-world dataset consists of gene expression data for 1027 female breast cancer (FBC) patients. An initial count of around 13k+ gene features were reduced down to 1146 using the disease-specific genomic analysis implemented by \citet{Nicolau2011}. To manually build the graph, the data was transformed to correlation distance and a lens function based on the standard deviation of each vector was applied. The parameters were set at 25 intervals with 80\% overlap. Clustering was performed using wards linkage with 3 clusters per interval. This Mapper graph revealed two flares, one highlighted by a small region of high survival at the tip composed of 2 nodes containing 3 tumour samples (Figure~\ref{fig:fbc_supplementary}). Following this, we reconstructed the graph using many different lens functions trying to find the same (or better) hotspot using our proposed algorithm. The lens functions we used were randomly generated as described in Section~\ref{sec:random_lenses}. For each scenario we ran the process (sampling of lens+Mapper+hotspot detection) 1000 times. As a criteria for hotspot detection we set the $\sigma_1$ at 5 for nodes, $\sigma_1$ at 10 for samples and $\epsilon$ at 0.15.
For each type of lens function, out of 1000 runs, we were able to find at least one graph with a hotspot region of increased survival. Our optimum result is shown in Figure~\ref{fig:linearsub}. The composition of each hotspot is described in Table~\ref{real-table}. 

\begin{table}
  \caption{%
  The results of hotspot detection on the TCGA FBC dataset. We describe the composition of each hotspot found for the various Mapper graphs constructed. "STD" refers to the manually created Mapper graph with a lens function built on standard deviation. Linear, linear subset, and non-linear subset refers to the three graphs constructed using the hotspot detection algorithm.}
  \label{real-table}
  \centering
  \begin{tabular}{lrrrr}
    \toprule
    
    \textbf{Lens}  & \textbf{STD} & \textbf{Linear} & \textbf{Linear subset} & \textbf{Non-linear subset} \\
    
    \midrule
    \textbf{n} & 3            & 51              & 21                     & 49                         \\
    \textbf{Survived (\%)}       & 100       & 98.04           & 100                 & 100                     \\
    \textbf{ER+ (\%)}            & 0            & 33.33            & 95.24                    & 85.71             \\
\textbf{ER- (\%)}               & 100          & 66.67            & 4.76                       & 14.29               \\
    \textbf{Basal (\%)}          & 100       & 41.18           & 0                   & 6.12                       \\
    \textbf{Her2 (\%)}           & 0         & 23.53           & 4.76                   & 4.08                       \\
    \textbf{LumA (\%)}           & 0      & 19.61           & 38.10                  & 65.31                      \\
    \textbf{Normal (\%)}         & 0         & 7.84            & 0.00                   & 16.33                      \\
    \textbf{LumB (\%)}           & 0         & 7.84            & 57.14                  & 8.16                     \\
    \bottomrule
  \end{tabular}
\end{table}

The manually created graph has a hotspot (Figure~\ref{fig:fbc_supplementary}) consisting of 3 ER- Basal patients (0.29\% of total cohort). Basal is a well-known subtype in breast cancer predicting poor survival prognosis and can easily be identified using standard hierarchical clustering \citep{Parker2009}. Our hotspot detection method was found to locate higher numbers of patients with increased survival in a single hotspot than the one presented in Figure~\ref{fig:fbc_supplementary}. When comparing the identity of patients, we found there was no overlap between the the manually created hotspot and the hotspots identified by our algorithm. Both the linear (Figure~\ref{fig:lin_supplementary}) and non-linear subset (Figure~\ref{fig:nonlibsub_supplementary}) lens functions detected 51 and 49 patients each respectively, but these were both determined to contain a mixture of ER status and breast cancer subtypes. The linear subset (Figure~\ref{fig:linearsub}) lens function contained a reduced number of samples, totalling 21 patients (2.05\% of the full cohort). Ninety-five percent of this patient hotspot (n=20) were found to have an ER+ status, and consisted solely of Her2, LumA, and LumB subtypes. This indicates an unusual cluster of patients similar to that found in \cite{Nicolau2011}.
\linespread{0.5}
\section{Conclusion}
In this paper we proposed a new method for hotspot detection on Mapper graph. This method could, for example, support biomedical analysis for precision medicine, where it is important to identify small groups of distinct patients that may exhibit varying behaviour. We demonstrated that the method worked well with three artificial Mapper graphs. Furthermore, for the TCGA real world dataset, the algorithm allowed us to construct lens functions, which led to graphs with hotspots. To further evaluate the algorithm, we will biologically validate the quality of the retrieved hotspots and investigate the presence of these hotspots in a secondary female breast cancer dataset. As a future direction for this work we will  explore the problem of overfitting while sampling from a space of lenses with the objective of hotspot detection.  

\bibliographystyle{natbib}
\bibliography{neurips_2020_ref}

\setcounter{section}{0}
\renewcommand{\thesection}{S\arabic{section}} 

\setcounter{figure}{0}
\renewcommand{\thefigure}{S\arabic{figure}}

\section{Supplementary Material}

\subsection{Definition of Mapper algorithm}
\label{sec:mapper}
In this section we introduce the Mapper algorithm originally introduced in~\cite{Singh2007}. For the convenience, the parameters of the algorithm are denoted with boldface. 

Let $X=\{x_1, \dots, x_k\}, x_i \in \mathbb{R}^n$ be the considered point cloud and and $A : X \rightarrow \mathbb{R}$ indicate the attribute of data (e.g. survival of patients). In addition let us consider a \textbf{\emph{lens} function} $f : X \rightarrow \mathbb{R}^n$. Typically $n=1$ are used. Let us take the range of $f(X)$ and cover it with \textbf{n} intervals that overlap on \textbf{k} percentage of their length. For instance, if $f(X) = [0,1]$, \textbf{n}=4 and \textbf{k}=20\%, the coverage will consist of the following four intervals of a length $\frac{5}{17}$: $[0,\frac{5}{17}], [\frac{4}{17},\frac{9}{17}], [\frac{8}{17},\frac{13}{17}]$ and $[\frac{12}{17},1]$. Now, take any interval $I$ from the coverage of $f(X)$ and consider $f^{-1}(I) = \{ x \in X | f(x) \in I \}$. We run a \textbf{clustering algorithm} of our choice at $f^{-1}(I)$ and obtain a collection of clusters $C^I_1,\ldots,C^I_m$. They will correspond to vertices of the Mapper graph. Given two intervals $I,J$ covering $f(X)$, two vertices $v_1,v_2$ corresponding to clusters $C^I_i$ and $C^J_k$ will be joined by an edge if and only if $C^I_i \cap C^J_k \neq \emptyset$.

The vertices and edges defined above constitute the \emph{Mapper graph} $G=(V,E)$. We now define a function $\hat{A} : V \rightarrow \mathbb{R}$ as:
	\begin{equation}
    \hat{A}(v) = \frac{1}{|v|}\sum_{x\in v}A(x)
	\end{equation}
where $A(x)$ is the value of the attribute.

\begin{figure}
    \centering
    \includegraphics[scale=0.5]{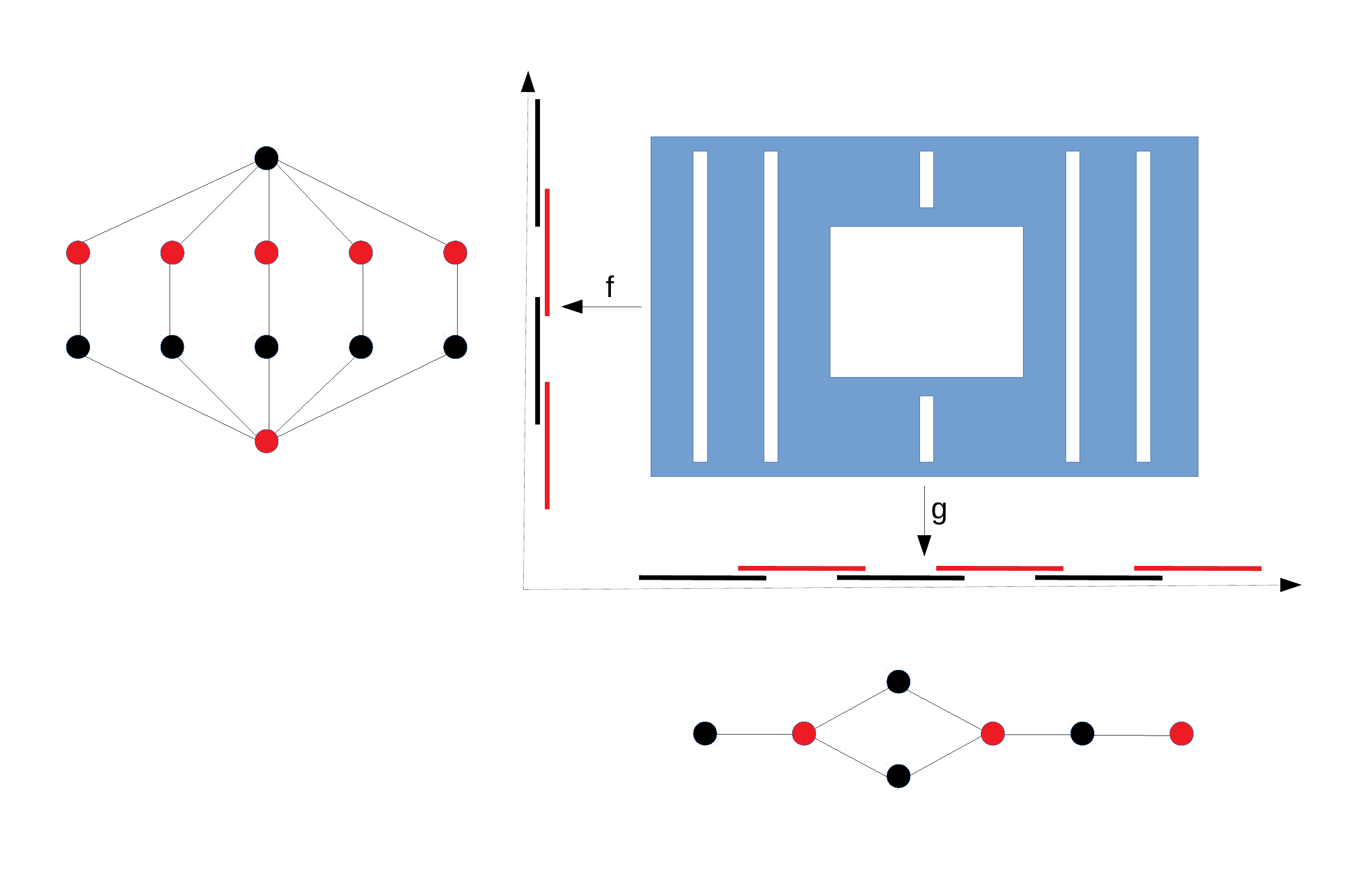}
    \caption{Example of two Mapper graph constructions. Consider the blue set as being very densely sampled with the point cloud $X$. In this figure we consider two possibilities of the lens functions: $f$ and $g$. The appropriate pullbacks and clustering therein provide the Mapper graphs on the left (for the function $f$) and in the bottom (for the function $g$).}
    \label{fig:different_lenses_mapper}
\end{figure}

A simple illustration of the presented construction can be found in the Figure~\ref{fig:different_lenses_mapper}.

Note that the example in the Figure~\ref{fig:different_lenses_mapper} highlights that different choices of lens functions may provide very different Mapper graphs. It is a natural phenomena, as the lens functions determine which details of the image are being neglected.

\subsection{Non uniqueness of components satisfying internal homogeneity criteria}
\label{sec:non_uniqueness_interlal_homogenity}
\begin{figure}[h!]
    \centering
    \includegraphics[scale=0.7]{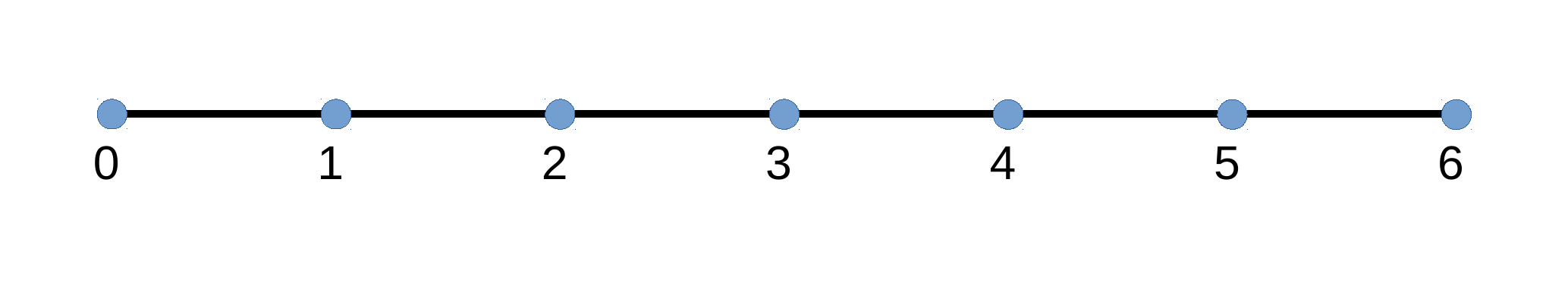}
    \caption{A graph with the value of the attribute function $A$ equal to the numbers below the vertices.}
    \label{fig:non_uniqueness_interlal_homogenity}
\end{figure}

The internal homogeneity criteria states that the value of the attribute function $A$ for any two vertices within a proposed region should not differ more than a chosen parameter $\tau$. Let us consider a graph together with the attribute function presented at Fig.~\ref{fig:non_uniqueness_interlal_homogenity}.

For a choice of $\tau = k$, for $k \in \{1,\ldots,6\}$, any connected subset of $k$ vertices satisfy the internal homogeneity condition. In such a case, an arbitrary subset is chosen in our algorithm. 

\subsection{Sampling from a space of lenses}
\label{sec:random_lenses}
Selection of lenses is often a nontrivial process. For data embedded to $\mathbb{R}^n$, the space of all possible functions is huge. Typically the considered lenses are constructed based on expert's opinion on the matter - on the knowledge that certain aspects of the data can be ignored (and therefore being put into fibers of the lenses). Still, in typical application, a number of lenses need to be considered before the one that is suitable for analysis is found.

The hotspot location procedures discussed in the Section~\ref{sec:hotspot_detection} opens a possibility of automatically sampling lens $f$ from the space of possible lenses and testing the quality of $f$. In this instance, the lenses giving cleaner hotspots will be considered better. 

In this instance, given the point cloud $X \subset \mathbb{R}^n$ and $(x_1,\ldots,x_n) \in X$ we consider random linear and quadratic functions of $x$. More precisely the considered linear functions are of the form:
\[ l(x) = \alpha_1 x_1 + \alpha_2 x_2 + \ldots + \alpha_n x_n  \]
where $\alpha_i$'s are randomly uniformly sampled from an interval $[-1,1]$. In some cases we consider most of $\alpha_i$'s equal to zero. Note that those lens functions are Lipschitz continuous and the presented process may be viewed as variable selection. 
We also consider quadratic functions of the form:
\[ q(x) = l(x) + \alpha_1 x_{r(\{1,\ldots,n\})} x_{r(\{1,\ldots,n\})} + \ldots + \alpha_1 x_{r(\{1,\ldots,n\})} x_{r(\{1,\ldots,n\})} \]
where $r(\{1,\ldots,n\})$ denotes a selection of random index among $\{1,\ldots,n\}$ and $\alpha_i$'s are, as above, random variables sampled from a uniform distribution.

\subsection{Graph Neighborhood}

\begin{figure}[h!]
\centering
\includegraphics[width=1\textwidth]{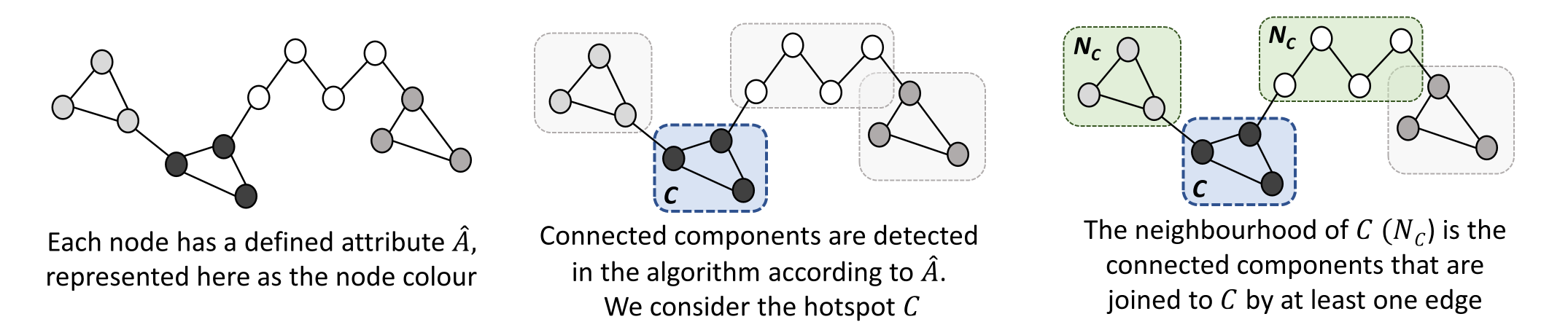} 
\caption{Overview of the approach to define the neighbourhood of a hotspot. An average attribute value $\hat{A}$ is assigned to each node. Using the cluster identification step of the hotspot detection algorithm, we identify four component groups in the Mapper graph. Considering whether a single group $C$ is a hotspot, the neighbourhood of $C$ is any connected component group that retains an edge with $C$. In this instance, two neighbourhood components are present. }
\label{fig:neighbourhood_sup}
\end{figure}

The definition of the Mapper graph neighbourhood as defined in Equation~\ref{eqn} of Section~\ref{PS} is summarised in Figure~\ref{fig:neighbourhood_sup}. 

\subsection{ Mapper graphs generated from hotspot detection algorithm}

\begin{figure}[h]
\centering
\includegraphics[width=0.6\textwidth]{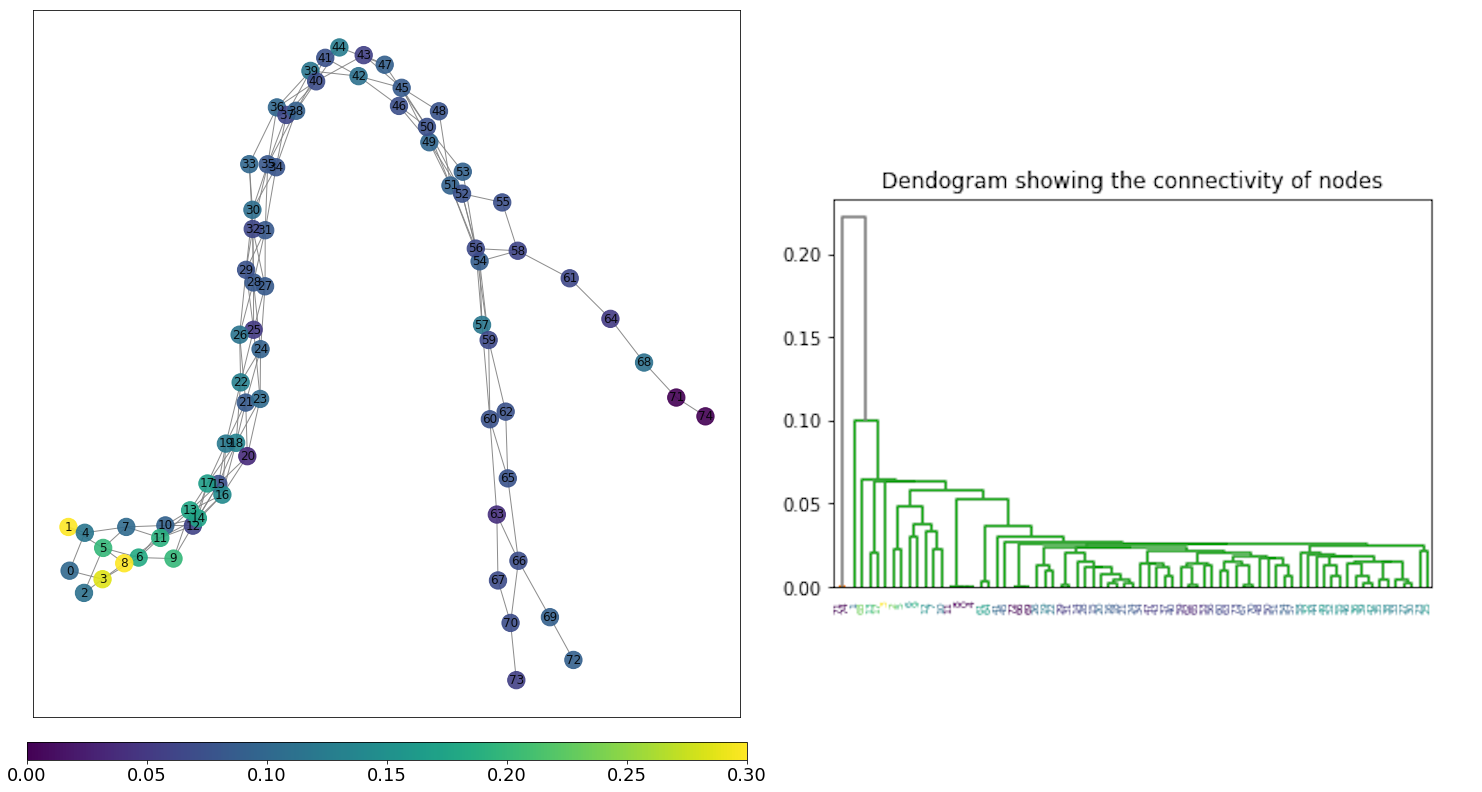} 
\caption{From left to right: Manually constructed Mapper graph created from the TCGA dataset, coloured by survival. The corresponding dendrogram of node connectivity based on edge weights. }
\label{fig:fbc_supplementary}
\end{figure}

\begin{figure}[h!]
\centering
\includegraphics[width=0.6\textwidth]{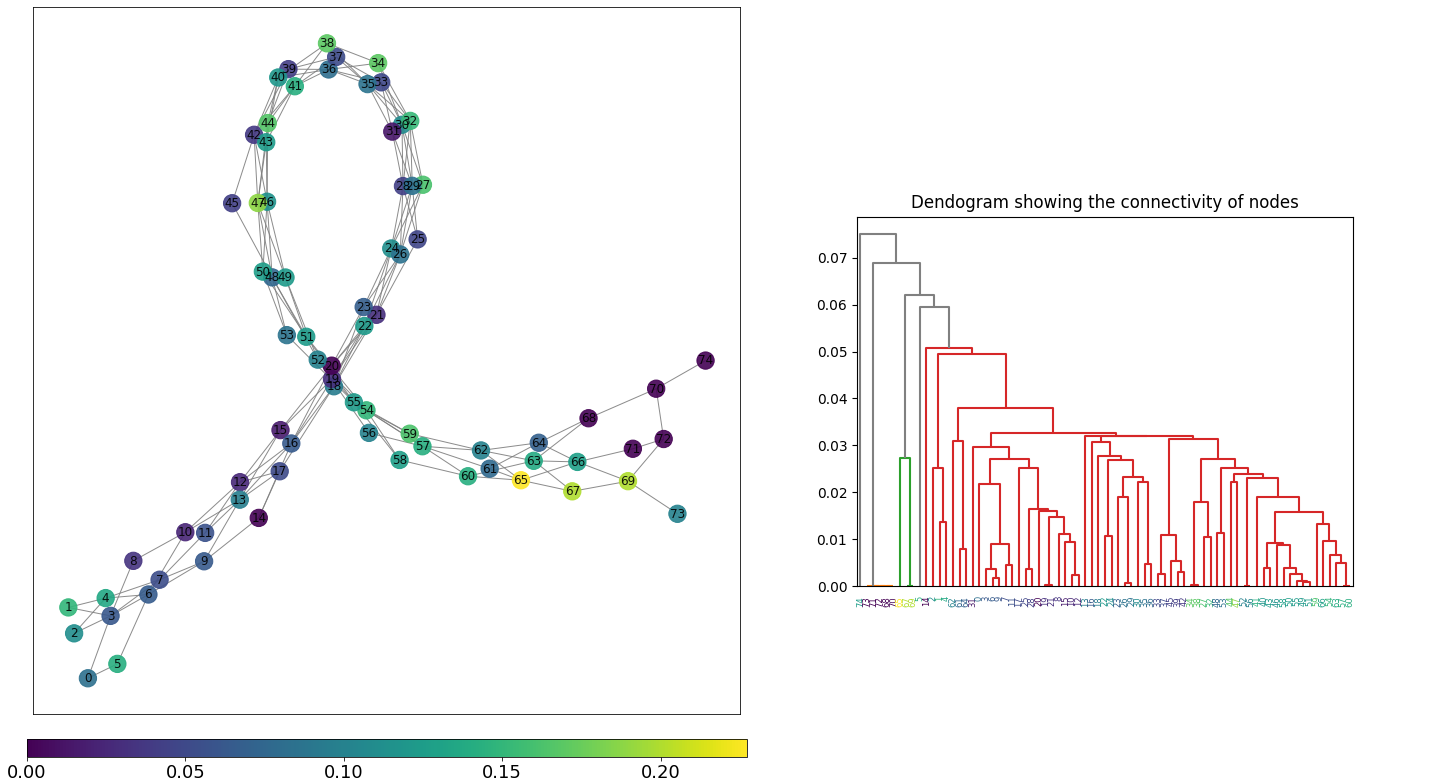} 
\caption{From left to right: The Mapper network graph built using a linear combination of features as the lens function, coloured by survival. The corresponding dendrogram of node connectivity based on edge weights.}
\label{fig:lin_supplementary}
\end{figure}

\begin{figure}[h!]
\centering
\includegraphics[width=0.6\textwidth]{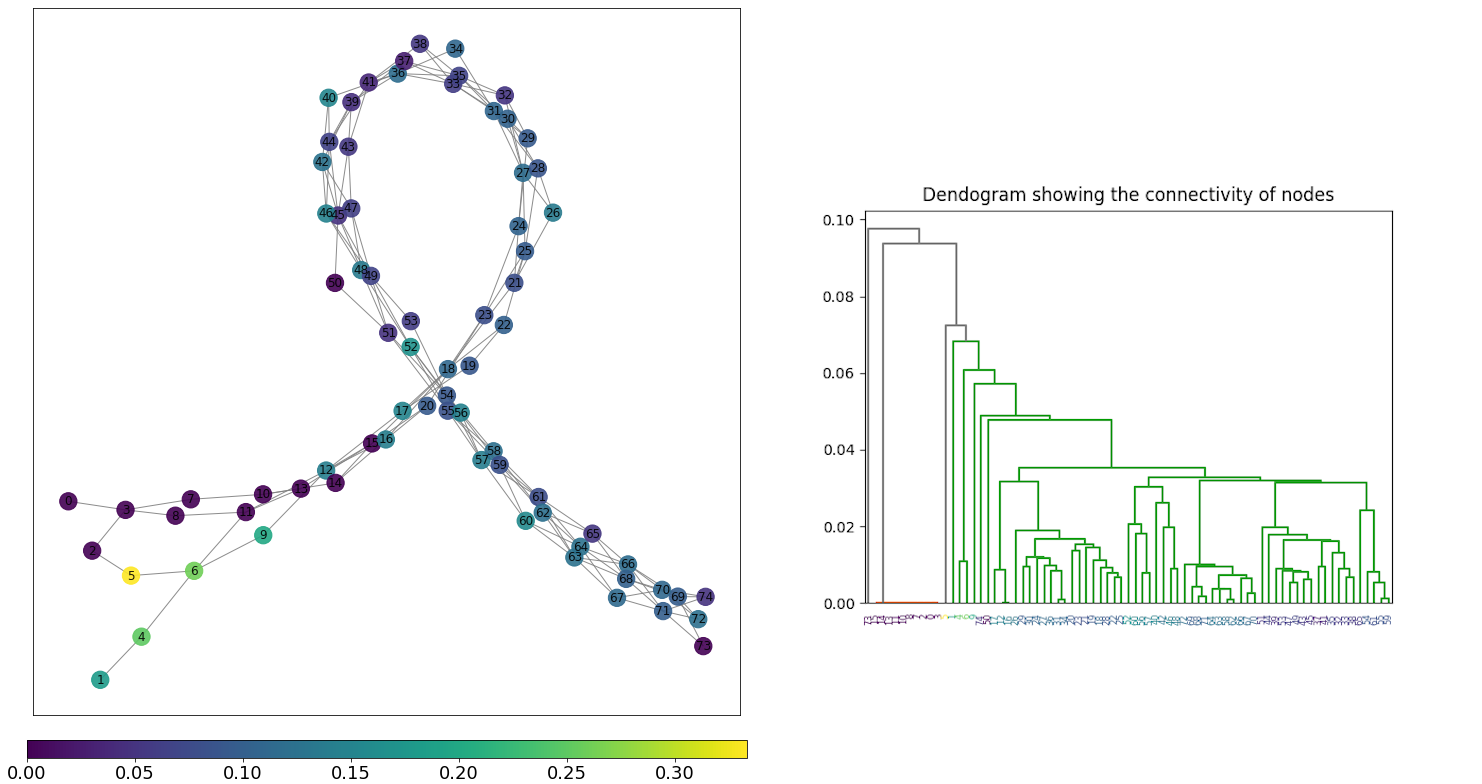} 
\caption{From left to right: The Mapper network graph built using a non-linear combination of subsetted features as the lens function, coloured by survival. The corresponding dendrogram of node connectivity based on edge weights.}
\label{fig:nonlibsub_supplementary}
\end{figure}

Please consider graphs on Figures~\ref{fig:fbc_supplementary}, \ref{fig:lin_supplementary} and \ref{fig:nonlibsub_supplementary} as examples of Mapper graphs that were obtained for random lenses sampled from the space of all lenses. 

\subsection{Stability with respect of selection of parameters}
Stability is a big issue of Mapper algorithm. Stability of the construction with respect to the coverage parameters (number of intervals and percentage of their overlap) was studied in~\cite{mapper_stability}. Yet clustering techniques used internally within the Mapper algorithm are inherently not stable, which makes the whole construction of the graph unstable in a general case. 

Consequently, classifications of sub-regions of the Mapper graph as a hotspot or not hotspot suffers from the same issue. 
It should be noted however that as a consequence some hotspots may be missed, but no \emph{false positive} answers will be obtained. This issue will be further elaborated in the full version of the paper. 

\section*{Acknowledgements}
This work was sponsored by a PhD studentship from the Northern Ireland Department for the Economy. PD acknowledge the support of Dioscuri program initiated by the Max Planck Society, jointly managed with the National Science Centre (Poland), and mutually funded by the Polish Ministry of Science and Higher Education and the German Federal Ministry of Education and Research.

\end{document}